\documentclass[a4paper,conference]{IEEEtran}
\usepackage{amsmath}
\usepackage{afterpage}
\usepackage{cite}
\allowdisplaybreaks[1]
\usepackage{graphicx}

\newcommand{\preeqn}{\vspace*{-0.03in}}
\newcommand{\preeqnt}{\vspace*{-0.02in}}
\newcommand{\preeqnm}{\vspace*{-0.02in}}

\newcommand{\posteqn}{}
\newcommand{\posteqnm}{}
\newcommand{\posteqnmm}{}
\newcommand{\posteqnl}{}

\usepackage[font=small,labelfont=bf,textfont=bf]{caption}

\renewcommand{\baselinestretch}{0.99}

\newcommand{\ul}{\underline}

\newcommand{\eqdef}{\buildrel \triangle \over =}
\newcommand{\xset}{{\cal R}^n}
\newcommand{\uset}{{\cal R}}
\newcommand{\yset}{{\cal R}^{2}}

\newcommand{\psiset}{{\cal R}^{n_\psi}}

\begin{document}
\IEEEoverridecommandlockouts

\title{
  Global Stabilization of Triangular Systems with Time-Delayed Dynamic Input Perturbations
}
\author{
  \IEEEauthorblockN{P. Krishnamurthy and F. Khorrami}
  \IEEEauthorblockA{Control/Robotics Research Laboratory (CRRL)\\
    Dept. of Electrical and Computer Engineering\\
    NYU Tandon School of Engineering\\
    Brooklyn, NY 11201, USA\\
    Email: \{prashanth.krishnamurthy,khorrami\}@nyu.edu}
  \thanks{\emph{To appear in 2017 IEEE International Carpathian Control Conference (ICCC).}}
}
\maketitle

\begin{abstract}
  A control design approach is developed for a general class of uncertain strict-feedback-like nonlinear systems with dynamic uncertain input nonlinearities with time delays.
  The system structure considered in this paper includes a nominal uncertain strict-feedback-like subsystem, the input signal to which is generated by an uncertain nonlinear input unmodeled dynamics that is driven by the entire system state (including unmeasured state variables) and is also allowed to depend on time delayed versions of the system state variable and control input signals.
  The system also includes
 additive uncertain nonlinear functions, coupled nonlinear appended dynamics, and uncertain dynamic input nonlinearities with time-varying uncertain time delays. 
 The proposed control design approach provides a globally stabilizing delay-independent robust adaptive output-feedback dynamic controller based on a dual dynamic high-gain scaling based structure.
\end{abstract}


\section{Introduction}
\noindent Consider the following class of uncertain nonlinear systems:
\preeqn\begin{align}
\dot x_i&=\phi_{(i,i+1)}(x_1)x_{i+1}+\phi_i(t,x)\,\, ,\,\,\, i=1,\ldots,n-1
\nonumber\\
\dot x_n&=  \mu(t,x,\psi,u)
\ \ \ ; \ \ \ \ 
\dot\psi = q_\psi(t,x,\psi,u,x_\Delta,\psi_\Delta,u_\Delta)
\nonumber\\
y &= [x_1,x_n]^T.
\label{eq:system}
\end{align}\posteqnmm
Here, the strict-feedback-like subsystem with state
$x=[x_1,\ldots,x_n]^T$ represents a nominal system and the subsystem with state
$\psi\in\psiset$ represents an appended input unmodeled dynamics.
$u\in\uset$ and $y\in\yset$ are the control input and measured output, respectively.
The subscript $\Delta$ is used to denote time delay, i.e., the notations $x_\Delta$, $\psi_\Delta$, and $u_\Delta$ refer to the time delayed versions of the signals $x$, $\psi$, and $u$, respectively, i.e.,  $x(t-\Delta)$, $\psi(t-\Delta)$, and $u(t-\Delta)$. Here, $\Delta$ is a (possibly time-varying) time delay\footnote{To simplify notation, the time argument is omitted when referring to signal values at time $t$, e.g., $u(t)$ is written simply as $u$; the time delayed signal values are written as $u_\Delta = u(t-\Delta)$, etc.}.
The functions $\phi_{(1,2)},\ldots,\phi_{(n-1,n)}$, are assumed to be known and continuous.
$\phi_1,\ldots,\phi_{n-1}$, $\mu$, and $q_\psi$ are uncertain continuous functions. 

In \eqref{eq:system}, the $x$ subsystem can be regarded as a
``nominal'' system, which if the value of $\mu$ could be directly specified by the controller, would be of the form 
$\dot x_i = \phi_{(i,i+1)}(x_1)x_{i+1} + \phi_i(t,x)\, ,\,
i=1,\ldots,n-1; \dot x_n =  u$.
In the actual system \eqref{eq:system}, the $x$ subsystem is driven by
a nonlinear uncertain function $\mu$, representing an input perturbation which involves $\psi$, which is the unmeasured state of the input unmodeled dynamics, as well as time delays versions of $x$, $\psi$, and $u$.
The control objective considered in this paper is to globally (i.e., starting from any initial condition) regulate the signals $x$ and $\psi$ in the system \eqref{eq:system} asymptotically to zero
under the various uncertainties described above and using measurement of the output $y$.

Control designs for various structures/classes of nonlinear dynamic systems including
 parametric and functional uncertainties, input nonlinearities and input unmodeled dynamics, time delays, etc., have been addressed
 in the literature (e.g., \cite{SOK84,KKK95,Isi99,Kha01,PW96,SJK97b,KK01a,KK01b,HLS07,RGLS08,CA09,Jan05,MB06,MML07,IR03,Pra03,KKC02,KK03b,KK03,Ito06,Lin09,KK09a,NK11} and references therein).
 Scaling based control designs for various types of triangular and non-triangular system structures have been addressed in \cite{KKC02,KK03a,KK04e,KK04f,KAP06,KK07f,KK07a,KK07b,KK07h,KK08a,KK10b,KK11b}.

 \IEEEpubidadjcol
 
 The application of the dynamic scaling technique to systems with uncertain input unmodeled dynamics such as in \eqref{eq:system} was considered in
 \cite{KK13a,KK14c,KK14b,KK15a}.
In \cite{KK13a,KK14c},
 nominal
 feedforward-like systems coupled with nonlinear input uncertainties were considered
 and a three-time-scale control design was developed that utilized two dynamic scaling parameters (one being essentially analogous to the scaling parameter that was utilized in our prior dual dynamic high-gain scaling based control designs \cite{KK04e} and the second scaling parameter being introduced specifically to handle the dynamic nonlinear input uncertainties).
 In \cite{KK15a}, the three-time scale (i.e., utilizing two dynamic scaling parameters) control design approach was extended to nominal strict-feedback-like systems coupled with nonlinear input uncertainties.
 In \cite{KK14b}, it was shown that (under certain structural conditions of the nominal system and nominal controller) a scaling based control redesign can be introduced to add robustness to input unmodeled dynamics to a general nonlinear system with a given nominal control design.
Here, we consider the general class \eqref{eq:system} of uncertain systems which includes time delays in the input unmodeled dynamics and we will show that the scaling-based design concept from \cite{KK13a,KK15a} can be applied to this uncertain system to provide global stabilization with robustness to uncertain input unmodeled dynamics. Specifically, while $q_\psi$ was required to be a function of $(t,x,\psi,u)$ at the current time $t$ in \cite{KK15a}, the control design here addresses the system structure shown in \eqref{eq:system} wherein $q_\psi$ involves time-delayed versions of both the state variables and the control input signals. The control design methodology developed in this paper is based on \cite{KK15a} and introduces refinements in the overall control design and the Lyapunov analysis to address the uncertain input unmodeled dynamics with time delays.

\renewcommand{\baselinestretch}{1.0}

\section{Notations and Assumptions}
\label{sec:assump}
\noindent{\bf Notations:} 
With $k$ being any integer,  the notation $I_k$ denotes an identity
matrix of dimension $k\times k$.
$|a|$ denotes Euclidean norm if $a$ is a column vector, absolute value if  
$a$ is a scalar, and Euclidean norm of the vector obtained by stacking all the columns of $a$ if $a$ is a (square or non-square) matrix. 
 Given any symmetric positive-definite matrix $P$, the notations
 $\lambda_{max}(P)$ and $\lambda_{min}(P)$ denote its maximum and minimum eigenvalues, respectively, of the matrix.
If $\alpha:[0,a)\rightarrow[0,\infty)$ is a strictly increasing continuous function with $\alpha(0)=0$, then
it is said to belong to class ${\cal K}$.
If, the class ${\cal K}$ definition holds with $a=\infty$ and
$\alpha(r)\rightarrow \infty$ as $r\rightarrow\infty$, then $\alpha$
is said to belong to class ${\cal K}_\infty$.

Here, we consider the output-feedback stabilization problem, i.e., only $y$ is assumed to be measured in the system \eqref{eq:system}.
It is assumed that the functions $\phi_i,\phi_{(i,i+1)},q_\psi,$ and $\mu$ appearing in the system dynamics satisfy sufficient conditions for local
existence and uniqueness of solution trajectories for the system (e.g., local Lipschitz conditions). 
The control objective in this paper is to make $x(t)$ and $\psi(t)$ in the system \eqref{eq:system}
asymptotically converge to zero
as $t\rightarrow\infty$
starting from any initial
conditions $x(0)\in{\cal R}^n$ and $\psi(0)\in S_\psi$. Here, $S_\psi$
is some known subset (possibly non-compact or unbounded) of ${\cal
  R}^{n_\psi}$. Also, $\overline S_\psi$ denotes the set of all
possible values of $\psi(t)$ over all time considering the set $S_\psi$ of possible initial values of $\psi$ and the
dynamics of the $\psi$ state variables.
Here, $\overline S_\psi$ could be simply ${\cal R}^{n_\psi}$ in general or could be a subset of ${\cal R}^{n_\psi}$.
The assumptions on the system considered here are given below.

\noindent{\bf Assumption A1:} (lower bound on magnitude of the {\em upper diagonal} terms $\phi_{(i,i+1)}$)
A positive constant $\sigma$ exists such that
$ \phi_{(i,i+1)}(x_1) \geq \sigma \,\, ,\,\, 1\leq i\leq n-1$ 
for all $x_1\in{\cal R}$.

\noindent{\bf Assumption A2:} (inequality bound on the uncertain functions $\phi_i$)
A known continuous function $\Gamma:{\cal R}\rightarrow{\cal R}^+$
and an unknown constant $\theta\geq 0$ exist
such that for all $t\geq 0$, $x\in\xset$, and $1\leq i\leq n$,
the following inequalities hold:
$|\phi_i(t,x)|
\leq
\Gamma(x_1)[\theta |x_1| + 
\sum_{j=2}^i|x_j|
]$.

\noindent{\bf Assumption A3:} (cascading dominance of upper diagonal terms)
The inequalities $\phi_{(i,i+1)}(x_1)  \geq  \overline \rho_i\phi_{(i-1,i)}(x_1)$
and
$\phi_{(i,i+1)}(x_1)  \leq  \ul\rho_i\phi_{(i-1,i)}(x_1)$ are satisfied
for all $x_1\in{\cal R}$ and $3\leq i \leq n-1$ with $\overline\rho_i$ and $\ul\rho_i$ being positive constants.

 \noindent{\bf Assumption A4:}  (assumptions on the uncertain nonlinear ``input perturbation'' function $\mu$)
 Known non-negative continuous
functions $\overline\mu(y,u)$,
$\overline\mu_1(y,u)$,
$\overline\mu_{1a}(x_1)$,
$\tilde\mu_1(y,u)$, and
$\overline\mu_2(y,u)$,
 non-negative (possibly
uncertain) functions
$\overline\mu_{1\psi}(\psi)$, $\tilde\mu_{1\psi}(\psi)$, and $\overline\mu_{2\psi}(\psi)$, an uncertain constant $\theta$,
and a known constant $\underline\mu$ exist such that for all time $t$,
 $x\in{\cal R}^n$, 
 $\psi\in \overline S_\psi$, and $u\in{\cal R}$:
(a) $\frac{\partial}{\partial u}\mu(t,x,\psi,u)
\geq \underline\mu>0$;
(b) $|\mu(t,x,\psi,u)|\leq \overline\mu(y,u)$;
(c) $|\frac{\partial \mu(t,x,\psi,u)}{\partial
  \psi}q_\psi(t,x,\psi,u,x_\Delta,\psi_\Delta,u_\Delta)|\leq
\overline\mu_1(y,u)
\sum_{k=0}^1\Big\{\overline\mu_{1a}(x_1(t-k\Delta))[\theta |x_1(t-k\Delta)| + \sum_{j=2}^n |x_j(t-k\Delta)|
]+ \overline\mu_{1\psi}(\psi(t-k\Delta))\Big\}$;
(d) $|\frac{\partial \mu(t,x,\psi,u)}{\partial t}|\leq \tilde\mu_1(y,u)[\theta
|x_1| + \sum_{j=2}^n |x_j|
+ \tilde\mu_{1\psi}(\psi)]$;
 (e) $|\frac{\partial \mu(t,x,\psi,u)}{\partial x}|\leq \overline\mu_2(y,u) + \overline\mu_{2\psi}(\psi)$.
 Also, for all $x\in\xset$, $t\in{\cal R}$, and $\psi\in\overline S_\psi$, we have $\lim_{|u|\rightarrow\infty}|\mu(t,x,\psi,u)|=\infty$.

\noindent{\bf Assumption A5:} (input-to-state stability -- ISS -- assumptions on the input unmodeled dynamics, i.e., the $\psi$ subsystem) 
A Lyapunov function $V_{\psi}:{\cal R}^{n_\psi}\rightarrow {\cal R}^+$ exists such that
for all $t\geq 0$,
$x\in{\cal R}^n$, $u\in{\cal R}$, and $\psi\in \overline S_\psi$, the following inequality holds:
$\frac{\partial V_\psi}{\partial\psi} q_\psi(t,x,\psi,u,x_\Delta,\psi_\Delta,u_\Delta) \leq
-\alpha_\psi(|\psi|)
+ \sum_{k=0}^1 \Big\{\Gamma_2(x_1(t-k\Delta)) [\theta x_1^2(t-k\Delta) +
\sum_{j=2}^n x_j^2(t-k\Delta)
+ \gamma_s(y(t-k\Delta),u(t-k\Delta))]\Big\}$,
with $\alpha_{\psi}$ being a known class ${\cal K}_\infty$ function,
$\gamma_s$ and $\Gamma_2$ being known non-negative functions, and
$\theta$ being an unknown non-negative constant.
Also,  $\gamma_s(y,u)\leq \overline\gamma_s(x_1)|\mu(t,x,\psi,u)|^2$ with $\overline\gamma_s$ being a known non-negative function.
Furthermore,
$\alpha_{\psi}(|\psi|) \geq \underline V_{\psi} V_\psi(\psi)$
and $[\overline\mu_{1\psi}^2(\psi) + \tilde\mu_{1\psi}^2(\psi) + \overline\mu_{2\psi}^2(\psi)] \leq \overline k_{\psi} \alpha_{\psi}(|\psi|)$
for all $\psi\in\psiset$
with $\underline
V_{\psi}$ and $\overline k_{\psi}$ being known positive constants 
and $\overline\mu_{1\psi}$, $\tilde\mu_{1\psi}$, and $\overline\mu_{2\psi}$ being functions of $\psi$ as in Assumption A4.

\noindent{\bf Assumption A6:} (assumptions on the time delay $\Delta$) The unknown time-varying time delay $\Delta$ is uniformly bounded in time and satisfies, for all time, the inequality given by $|\dot\Delta|\leq \overline\Delta < 1$ with $\dot\Delta$ denoting $\frac{d}{dt}\Delta$ and with $\overline\Delta$ being a known positive constant.

\noindent{\bf Theorem 1:}
Under the Assumptions A1-A6, 
a positive constant
$a_\theta$
and continuous functions
$g_i,i=2,\ldots,n$,
$k_i,i=2,\ldots,n$,
and nonnegative continuous functions
$\vartheta_1$,
$\lambda$, 
$R$, 
$\Omega$, 
$R_u$, 
$\Omega_u$,
 and 
$\overline Q_\theta$
can be found such that the dynamic output-feedback controller given below (with $[\hat x_2,\ldots,\hat x_n,\zeta,r,r_u,\hat\theta]^T$ comprising the state of the designed dynamic controller)
\preeqn\begin{align}
&\left.\begin{array}{rcl}
\dot{\hat x}_i&\!\!=&\!\!\phi_{(i,i+1)}(x_1)[\hat x_{i+1}
+r^i f_{i+1}(x_1)]
\\&\!\!&\!\!
-r^{i-1} g_i(x_1)[\hat x_2+r  f_2(x_1)]\\&\!\!&\!\!
-(i-1)\dot r r^{i-2} f_i(x_1)\,\,\, ,\,\,\, 2\leq i\leq n-1
\\
\dot{\hat x}_n&\!\!=&\!\!\tilde u
-r^{n-1} g_n(x_1)[\hat x_2+r  f_2(x_1)]
\\&\!\!&\!\!
-(n-1)\dot r r^{n-2} f_n(x_1)
\\
\hat x &\!\!=&\!\! [\hat x_2,\ldots,\hat x_n]^T
\\
f_i(x_1) &\!\!=&\!\! \int_0^{x_1}\frac{g_i(\pi)}{\phi_{(1,2)}(\pi)}d\pi \ ,\ i=2,\ldots,n
\end{array}\right\}
\label{obsdyn}
\\
&\!\!\!\!\left.\begin{array}{rcl}
\varpi_2 &\!\!\!\!=&\!\!\!\! \frac{\hat x_2+rf_2(x_1)+\vartheta (x_1,\hat\theta)}{r}
\\
\varpi_i &\!\!\!\!=&\!\!\!\! \frac{\hat x_i+r^{i-1}f_i(x_1)}{r^{i-1}}
\,\,\, ;\,\,\,
i=3,\ldots,n
       \end{array}\!\!\!\right\} \ ;  \ 
          \varpi = [\varpi_2,\ldots,\varpi_{n}]^T
\label{xiidefn}
\end{align}\posteqnm
       \vspace*{-0.2in}
\preeqn\begin{align}
\vartheta(x_1,\hat\theta) = \hat\theta x_1 \vartheta_1(x_1)&\qquad\qquad\qquad\qquad\qquad\qquad\qquad
\label{varthetadefn}
       \end{align}\posteqnm
       \vspace*{-0.2in}
\preeqn\begin{align}
         \tilde u&= -r^n K(x_1)\varpi
                   \ \ ; \ 
                   K(x_1) \!\eqdef\! [k_2(x_1),k_3(x_1),\ldots,k_{n}(x_1)]
\label{tilde_u_defn}
\\
\dot \zeta &= \dot r_u x_n +
r_u \tilde u
\,\,\,\, ;\,\,\,
        u 
        =
        \zeta - r_u x_n
	\label{udefn1}
\\
 \dot r &=
\lambda\Big(R\Big(x_1,\hat\theta,\dot{\hat\theta})-r\Big)
\Omega(x_1,\hat\theta,\dot{\hat\theta},r)
\,\,\, ;\,\,\,
r(0) \geq 1
\label{rdot}
\\
\dot r_u &=
\lambda(R_u(y,u,r,\dot r,\hat\theta,\dot{\hat\theta},\varpi)-r_u)
\nonumber\\&\quad
\times
\Omega_u(y,u,r,\dot r,\hat\theta,\dot{\hat\theta},\varpi,r_u)
\,\,\, ;\,\,\,
r_u(0) \geq 1
\label{r2dot} 
\\
  \dot{\hat\theta} &=
\overline Q_\theta(x_1) 
\,\,\, ;\,\,\,
\hat\theta(0) \geq a_{\theta}
\label{thetahdot}
\end{align}\posteqnmm
when put in closed loop with the system with dynamics shown in \eqref{eq:system} guarantees that
for any initial conditions $(x(0), \psi(0), \hat x(0)$, $\zeta(0),r(0),r_u(0),\hat\theta(0))\in{\cal R}^n \times S_\psi \times {\cal R}^{n-1} \times {\cal R}\times [1,\infty)\times [1,\infty)\times [a_\theta,\infty)$, all the closed-loop signals ($x_1,\ldots,x_n$, $\psi$, $\hat x_2,\ldots,\hat x_n$, $\zeta$, $r$, $r_u$, $\hat\theta$, $u$) are uniformly bounded over the time interval $[0,\infty)$ and, furthermore, the signals 
$x_1,\ldots,x_n,\psi,\hat x_2,\ldots, \hat x_n$ converge to zero asymptotically as the time $t$ goes to $\infty$.

\noindent{\bf Remark 1:} Comparing assumptions A1 through A6 given above with the corresponding assumptions in
our earlier work \cite{KK15a}, it is seen that
A1 through A5 are essentially analogous to the
assumptions considered before.
The additional element in the system structure here is the presence of time delays in the input unmodeled dynamics. The required assumption on this additional element is addressed by Assumption A6.
Assumption A6 is equivalent to the standard assumption utilized in the
literature on time delayed systems (e.g., \cite{KK14c}) that essentially requires that the
time delay value does not change faster than ``real-time'' (i.e.,
$|\dot\Delta| < dt/dt = 1$).
The proposed dynamic controller design approach
can be applied to systems that also have time
delays in other parts of the system (e.g., in the $\phi_i$ terms) and multiple time delay values (e.g.,
state and input time
delays in the nominal system dynamics,
multiple possible time delay values instead of a single $\Delta$, etc.) by appropriately adding additional terms in the overall system Lyapunov function. However,
these additional possible time delays are not considered here so as
to focus on the basic control design approach to
handle the time delay $\Delta$ in the input unmodeled dynamics.
$\diamond$

\noindent{\bf Remark 2:}
The proposed control design given in Theorem 1 comprises
of
a reduced-order observer
(with state vector $[\hat x_2,\ldots, \hat x_n]^T$)
with dynamics in \eqref{obsdyn}, the nominal control law \eqref{tilde_u_defn}, a dynamic state extension $\zeta$ with dynamics as in \eqref{udefn1}, dynamic scaling parameters $r$ and $r_u$ with dynamics as in \eqref{rdot} and \eqref{r2dot}, respectively, and an adaptation parameter $\hat\theta$ with dynamics as in \eqref{thetahdot}.
Here, the dynamic scaling parameters $r$ and $r_u$ are initialized with $r(0) \geq 1$
and $r_u(0) \geq 1$ and the dynamic adaptation parameter $\hat\theta$ is initialized with
$\hat\theta(0) \geq a_\theta > 0$.
From the dynamics of these controller state variables,
we see that $\dot r(t)$, $\dot r_u(t)$, and $\dot{\hat\theta}(t)$ are non-negative
 at all time $t$.%
It is noteworthy that the overall dynamic controller structure in Theorem~1 is essentially as in \cite{KK15a} and it is seen in Section~\ref{sec:redesign} that the time delays in the input unmodeled dynamics are handled through changes in the designs of the functions $R$, $\Omega$, $R_u$, $\Omega_u$, etc., while retaining the overall controller structure. This is indeed illustrative of the flexibility and generality of the dynamic scaling-based controller design approach that enables (as noted in, for example \cite{KK07a,KK14b}) control designs for a wide range of classes of systems and uncertainty structures within a unified framework.   
$\diamond$

\vspace*{0.03in}
\noindent{\bf Remark 3:}
Consider the following system with output
$y = [x_1,x_4]^T$:
\preeqn\begin{eqnarray}
\dot x_1 &\!\!\!\! =&\!\!\!\!  (1+x_1^2)x_2  + \theta_1 x_1^2\cos(x_3)
\nonumber\\
 \dot x_2 &\!\!\!\! =&\!\!\!\! (1+x_1+x_1^2)x_3  + \theta_2 x_1^3 \cos(x_2) +x_1^2 x_2
\nonumber\\
 \dot x_3 &\!\!\!\! =&\!\!\!\!  (1+2x_1^2)x_4  + \theta_3 x_1^2 + b_1 (x_1^3 x_3 + x_1 x_4 )
\nonumber\\
\dot x_4 &\!\!\!\! =&\!\!\!\! \mu(t,x,\psi,u)
\nonumber\\
&\!\!\!\!=&\!\!\!\! a_1(2+\cos(\psi_1))u \!+\! a_1\sin(t)x_1u^2
\nonumber\\&\!\!\!\!&\!\!\!\!
+ (a_1+a_2) (2+\sin(x_3)\sin(\psi_2)+x_1^2+x_4^2)u^3
\nonumber\\
\dot\psi_1 &\!\!\!\! =&\!\!\!\!  \psi_2 - \psi_1 + b_2 x_{1,\Delta} \sin(x_3) + b_3 \cos(t) x_1 x_{1,\Delta}^2 
\nonumber\\&\!\!\!\!&\!\!\!\!
+ b_4 x_{1,\Delta}^2 x_{3,\Delta} \cos(x_{4,\Delta}) + b_5 u \cos(\psi_{1,\Delta})x_1
\nonumber\\
         \dot\psi_2 &\!\!\!\! =&\!\!\!\!  -2\psi_2 +\psi_2\cos(\psi_{1,\Delta}\psi_{2,\Delta}) + b_6 x_1 x_{2,\Delta} \sin(x_4)
                                 \nonumber\\&&
                                 + b_7 x_1 u\cos(u_\Delta)
\label{exsys}
\end{eqnarray}\posteqnmm
where $x_{i,\Delta}$ and $\psi_{i,\Delta}$ denote $x_i(t-\Delta)$ and $\psi_i(t-\Delta)$, respectively.
$\theta_i,i=1,2,3$, $b_i,i=1,\ldots,7$, and $a_i,i=1,2$, are unknown constants; we assume that upper and lower bounds are known for $a_i,i=1,2$, as $\overline a_i$ and $\underline a_i > 0$, respectively.
and upper bounds $\overline b_i,i=1,\ldots,7$ are known for $|b_i|$. 
Also, an upper bound $\overline\Delta < 1$ on $|\dot\Delta|$ is considered as known.
  This sixth-order system can be seen to satisfy the Assumptions A1-A6. The value of $\sigma$ in Assumption A1 can be picked as $3/4$.
  The function $\Gamma$ in
Assumption A2 can be chosen to be $\Gamma(x_1)=\max(1,\overline b_1)|x_1|+x_1^2+\overline b_1 x_1^3$ and the uncertain parameter $\theta$ can be defined as $\theta_a= \max(1,|\theta_1|,|\theta_2|,|\theta_3|)$.
With $\overline\rho_3=0.8$, and $\underline\rho_3=4$, we see that Assumption A3 is also satisfied.
Also, with 
$\underline\mu = \frac{2}{3}\underline a_1$, 
$\overline\mu(y,u) =
3\overline a_1 |u| \!+\! \overline a_1|x_1|u^2 \!+\! (\overline a_1+\overline a_2) (3+x_1^2+x_4^2)|u|^3$,
$\overline\mu_1(y,u) =
[\overline a_1 |u| + (\overline a_1+\overline a_2)|u|^3][1+\overline b_2  + \overline b_3 |x_1|  + \overline b_4   + \overline b_5  |u|
+ \overline b_6 |x_1|  + \overline b_7 |u|]$,
$\overline\mu_{1a}(x_1)=1+|x_1|+x_1^2$,
$\overline\mu_{1\psi}(\psi) = |\psi_1|+3|\psi_2|$,
$\tilde\mu_1(y,u) = \overline a_1 u^2$, $\tilde\mu_{1\psi}(\psi) = 0$,
$\overline\mu_2(y,u) = \overline a_1 u^2 + (\overline a_1+\overline a_2)(1+2|x_1|+2|x_4|) |u|^3$, 
$\overline\mu_{2\psi}(\psi) = 0$, and with $\theta$ appearing in Assumption A4 defined as $\theta_b = 1$, it is seen that Assumption A4 is satisfied.
 Defining $V_\psi = \frac{1}{2}(\psi_1^2+\psi_2^2)$, Assumption A5 holds with
$\alpha_{\psi} (|\psi|) = \frac{1}{4}|\psi|^2$,
$\theta$ given as $\theta_b = 1$,
$\Gamma_2(x_1)=4(\overline b_2^2 + \overline b_3^4x_1^2 + \overline b_3^4x_1^6 + \overline b_4^2x_1^4 + \overline b_5^2) + 2(\overline b_6^4x_1^2  + \overline b_7^2)$,
$\gamma_s(y,u)=(1+x_1^2)^2 u^2$,
$\overline\gamma_s(x_1) = [4(1+x_1^2)/(3\underline a_1)]^2$,
$\overline k_{\psi} = 36$, and $\underline V_{\psi} = \frac{1}{2}$.
The overall uncertain parameter $\theta$  appearing in Assumptions~A2 and A5 can be defined as $\theta = \max(\theta_a,\theta_b)$.
Therefore, as described above, the example system \eqref{exsys} satisfies the Assumptions A1-A6 and the proposed control design methodology can be applied to this example system.
$\diamond$

\section{Control Design and Stability Analysis
}
\label{sec:redesign}
\noindent {\bf Choosing $g_i,k_i,i=2,\ldots,n$:}
Define
$D_o= D_c = \mbox{diag}(1,2,\ldots,n-1)$
              and
              $\tilde D_o = \tilde D_c = D_o \!- \frac{1}{2}I_{n-1}$.
Let $C$ be the $1\times (n-1)$ matrix given as $[1,0\ldots,0]$.
Let $A_o(x_1)$ and $A_c(x_1)$ be the square matrices of dimension $(n-1)\times (n-1)$ with $(i,j)^{th}$ elements defined as
  $A_{o_{(i,i+1)}} = \phi_{(i+1,i+2)} \,,\, i=1,\ldots,n-2$,
  $A_{o_{(i,1)}} = -g_{i+1} \,,\, i=1,\ldots,n-1$,
   $A_{c_{(i,i+1)}} = \phi_{(i+1,i+2)}\, , \, i=1,\ldots,n-2$ and
 $A_{c_{(n,j)}} = -k_{j+1}
 \, ,\, j=1,\ldots,n-1$, and with all other elements being zero.
 Given Assumptions A1 and A3, we know  (\cite{KK04f,KK04g,KK06g}) that
 constant $(n-1)\times (n-1)$ matrices $P_o>0$ and $P_c>0$,  positive constants $\nu_{o}$, $\tilde\nu_{o}$ $\underline\nu_{o}$,
 $\overline\nu_{o}$, $\nu_{c}$, $\underline\nu_{c}$,
and $\overline\nu_{c}$, and functions $g_2,\ldots,g_n$,$k_2,\ldots,k_n$ can be constructed such that the
following two pairs of coupled Lyapunov inequalities are satisfied for all $x_1\!\in\!{\cal R}$:  
\preeqn\begin{eqnarray}
&P_o A_o(x_1)+A_o^T(x_1) P_o\leq-\nu_o I_{n-1} -\tilde\nu_o \phi_{(2,3)}(x_1) C^T C&\nonumber\\
&\underline \nu_o I_{n-1} \leq P_o\tilde D_o+\tilde D_o P_o\leq\overline \nu_o I_{n-1} &
\label{coupledlyapeqns1}
       \end{eqnarray}\posteqnmm
       \vspace*{-0.2in}
       \preeqn\begin{eqnarray}
&P_c A_c(x_1)+A_c^T(x_1) P_c\leq-\nu_c \phi_{(2,3)}(x_1)I_{n-1}&\nonumber\\
&\underline \nu_c I_{n-1} \leq P_c\tilde D_c+\tilde D_c
P_c\leq\overline \nu_c I_{n-1}. &
\label{coupledlyapeqns2}
              \end{eqnarray}\posteqnmm
              The first of these pairs \eqref{coupledlyapeqns1} can be considered the {\em observer-context} coupled Lyapunov inequalities that relates to the choice of the functions ({\em observer gain functions}) $g_2,\ldots,g_n$ while the second of these pairs \eqref{coupledlyapeqns2} can be considered the {\em controller-context} coupled Lyapunov inequalities that relates to the choice of the functions ({\em controller gain functions}) $k_2,\ldots,k_n$. 
 Also, the functions $g_2,\ldots,g_n$ can be chosen (\cite{KK04f,KK04g,KK06g}) 
such that  
$\sqrt{\sum_{i=2}^n g_i^2(x_1)} \leq \overline g\phi_{(2,3)}(x_1)$ for all $x_1\in {\cal R}$
with $\overline g$ being a positive constant.

\vspace*{0.03in}
\noindent{\bf Scaled observer errors $\epsilon_i$ and their dynamics:}
As noted in Remark 2, the dynamics \eqref{obsdyn} can be viewed as a  reduced-order observer with
state variables $\hat x=[\hat x_2,\ldots,\hat x_n]^T$.
Define the
observer error vector $e = [e_2,\ldots,e_n]^T$ and the scaled
observer error vector $\epsilon =
[\epsilon_2,\ldots,\epsilon_n]^T$ as
$e_i = \hat x_i+r^{i-1}f_i(x_1)-x_i\,\, ,\,\,\, i=2,\ldots,n$ and
$\epsilon_i=\frac{e_i}{r^{i-1}}\,\,\, ,\,\,\, i=2,\ldots,n$.
Then, we have the dynamics
\preeqn\begin{align}
         \dot\epsilon &= r  A_o \epsilon-\frac{\dot r}{r}D_o\epsilon\!+\!\overline\Phi \!-\! B \frac{u_d}{r^{n-1}}
                        \ \ ; \ \ 
                        u_d \eqdef \mu(t,x,\psi,u)\!-\!\tilde u
\label{obs_scaled_error_dynamics}
\end{align}\posteqnm
where
$B\in{\cal R}^{n-1} = [0,\ldots,0,1]^T$, and $\overline\Phi=[\overline\Phi_2,\ldots,\overline\Phi_n]^T$ with
$\overline\Phi_i=-\frac{\phi_i}{r^{i-1}}+g_i\frac{\phi_1}{\phi_{(1,2)}}$.
The quantities $\varpi_i$ in \eqref{xiidefn} can be viewed as
scaled observer estimates with dynamics given by
  $\dot\varpi = r A_c \varpi 
-\frac{\dot r}{r}D_c\varpi-r G
\epsilon_2+G\frac{\phi_1}{\phi_{(1,2)}} +B_1\Psi$
where
$G=[g_2,\ldots,g_n]^T$,
$B_1\in{\cal R}^{n-1}$ is given as  $[1,0,\ldots,0]^T$, and
$\Psi = 
    \frac{\partial\vartheta}{\partial x_1}\phi_{(1,2)}(x_1)(\varpi_2-\epsilon_2)
    + \frac{1}{r}\frac{\partial\vartheta}{\partial x_1}[\phi_1-\phi_{(1,2)}\vartheta]
    +\frac{1}{r}\frac{\partial\vartheta}{\partial\hat\theta}\dot{\hat\theta}$.
    From
    \eqref{varthetadefn}, we have
$\frac{\partial\vartheta}{\partial\hat\theta} = x_1\vartheta_1(x_1)$ and
$\frac{\partial\vartheta}{\partial x_1} = \hat\theta(\vartheta_1(x_1)+x_1\vartheta_1'(x_1))$ where the notation $\vartheta_1'(x_1)$ is used to represent $\frac{d\vartheta_1(\pi)}{d\pi}\Big|_{\pi=x_1}$.

\noindent{\bf Observer Lyapunov function $V_o$, controller Lyapunov function $V_c$, and composite Lyapunov function $V_x$:}
The observer and controller Lyapunov functions are defined as
$V_o = r\epsilon^T P_o \epsilon$ and
$V_c = r\varpi^T P_c \varpi+\frac{1}{2}\left(1+\frac{1}{r}\right)x_1^2$
and a composite Lyapunov function is defined as the linear combination
$V_x = c V_o + V_c$
where 
$c > 0$ is any constant such that
$c > [4 \lambda_{max}^2(P_c)\overline g^2/(\tilde\nu_o\nu_c)]$.
Using \eqref{coupledlyapeqns1} and \eqref{coupledlyapeqns2}, it can be shown that
\preeqn\begin{align}
         \dot V_x
  &\leq
-r^2\frac{c\nu_o}{4}|\epsilon|^2
\!-\!r^2\frac{\nu_c}{4}\phi_{(2,3)}(x_1)|\varpi|^2
\!-\!\dot r c\ul\nu_o|\epsilon|^2
\!-\!\dot r \ul\nu_c|\varpi|^2
\nonumber\\&\quad
-x_1\phi_{(1,2)}\vartheta 
- \frac{1}{2}\frac{\dot r}{r^2} x_1^2 
  - \frac{2c r\epsilon^T P_o B u_d}{r^{n-1}} 
\nonumber\\&\quad
  + rw_1(x_1,\hat\theta,\dot{\hat\theta})[|\epsilon|^2\!+\!|\varpi|^2] 
\!+\!q_1(x_1)x_1^2 \!+\! \theta^* q_2(x_1) x_1^2
\label{vxdot}
\end{align}\posteqnm
where $\theta^* = \theta + \theta^2$,
$q_1(x_1) = \frac{4}{c\nu_o}\phi_{(1,2)}^2(x_1) +
    \frac{8}{\nu_c \phi_{(2,3)}(x_1)}\phi_{(1,2)}^2(x_1)
    + \lambda_{max}^2(P_c)
    + 2$,
    $q_2(x_1) = 2\Gamma(x_1) + 1
    + \frac{8\lambda_{max}^2(P_c)\overline g^2
      \phi_{(2,3)}(x_1)}{\nu_c}\frac{\Gamma^2(x_1)}{\phi_{(1,2)}^2(x_1)} 
    +\frac{8n}{c\nu_o}\lambda_{max}^2(P_o)\Gamma^2(x_1)\Bigg[1+\overline
      g^2 \phi_{(2,3)}^2(x_1)/\phi_{(1,2)}^2(x_1)\Bigg]$,
    and\\
    $w_1(x_1,\hat\theta,\dot{\hat\theta}) = 3\lambda_{max}(P_c)
    \Big|\frac{\partial\vartheta(x_1,\hat\theta)}{\partial x_1}\Big|\phi_{(1,2)}(x_1) +
    \vartheta_1^2(x_1)\dot{\hat\theta}^2
    +
    \lambda_{max}^2(P_c)\hat\theta^2(\vartheta_1(x_1)+\vartheta_1'(x_1)x_1)^2
(\Gamma^2(x_1)
    + \phi_{(1,2)}^2(x_1)\vartheta_1^2(x_1))
    + 3\lambda_{max}(P_o)(n+n^2)\Gamma(x_1)
    + n\lambda_{max}^2{P_o}\hat\theta^2\vartheta_1^2(x_1)$.
    
\vspace*{0.03in}
\noindent{\bf Remark 4:}
The nominal $x$ subsystem can be seen to be globally stabilized if the virtual ``control input'' entering into the subsystem could be made to be $\tilde u$ instead of $\mu(t,x,\psi,u)$. Hence,
$u_d =
[\mu(t,x,\psi,u)-\tilde u]$ can be viewed as a ``mismatch'' term due to the fact that the actual control signal entering into the $x$ subsystem is  $\mu(t,x,\psi,u)$ instead of the desired virtual control input $\tilde u$.
It will be seen in the analysis below that
the control ``redesign'' given by the dynamic state extension $\zeta$ and control input ($u$) definition in
  \eqref{udefn1} along with the designs of dynamics of the scaling parameters $r$ and $r_u$ will make $\mu(t,x,\psi,u)$ track the desired/nominal control input signal
  $\tilde u$, thus making the actual system \eqref{eq:system} with the input unmodeled dynamics  globally stabilized with asymptotic convergence of $x$ and $\psi$.
  $\diamond$ 

\noindent{\bf Mismatch term $u_d =
  [\mu(t,x,\psi,u)-\tilde u]$ and Lyapunov function component $V_u$:}
We see from \eqref{udefn1} that:
\preeqn\begin{align}
	\dot u &= \dot \zeta - \dot r_u x_n - r_u \dot x_n
        =
     -r_u [\mu(t,x,\psi,u)-\tilde u] 
 =  -r_u  u_d.
	\label{udefn2}
\end{align}\posteqnm
Hence,
\preeqn\begin{align}
  \dot u_d &=
-r_u  \frac{\partial\mu(t,x,\psi,u)}{\partial u}u_d
+ \frac{\partial\mu(t,x,\psi,u)}{\partial x_n}u_d
+\chi_1
\label{uddot}
\end{align}\posteqnm
where
\preeqn\begin{align}
  \chi_1 &= 
  \frac{\partial\mu(t,x,\psi,u)}{\partial x_1}[\phi_{(1,2)}(x_1) x_2 + \phi_1(t,x)]
           \nonumber\\&
                          +
  \frac{\partial\mu(t,x,\psi,u)}{\partial x_n}\tilde u
               \!+\! n r^{n-1}\dot r K(x_1)\varpi
           \nonumber\\&
               +\!\!\sum_{i=2}^{n-1}\!\!\frac{\partial\mu(t,x,\psi,u)}{\partial x_i}[\phi_{(i,i+1)}(x_1)x_{i+1} \!+\! \phi_i(t,x)]
  \nonumber\\&
               + \frac{\partial\mu(t,x,\psi,u)}{\partial\psi} q_\psi(t,x,\psi,u,x_\Delta,\psi_\Delta,u_\Delta)
  \nonumber\\&
               + \frac{\partial\mu(t,x,\psi,u)}{\partial t} 
+ r^n \frac{\partial K}{\partial x_1}\varpi [\phi_{(1,2)}(x_1) x_2 + \phi_1(t,x)]
               \nonumber\\&
+ r^n K(x_1)\dot\varpi.
\label{chi1defn}
\end{align}\posteqnm
Here, $\frac{\partial\mu(t,x,\psi,u)}{\partial x_i}$
denotes $\frac{\partial\mu(t,x,\psi,u)}{\partial
  x}\varrho_i$ where $\varrho_i$ is a vector of dimension $n\times 1$ having 1
 as its $i^{th}$ element and having 0's everywhere else.
The quantity $\chi_1$ can be shown to satisfy the magnitude bound
    $|\chi_1| \leq 
    \beta_1(y,u,r,\dot r,\hat\theta)[|\epsilon| + |\varpi|]
        + \beta_2(y,u,r,\hat\theta,\dot{\hat\theta}) |x_1| 
        + \theta \beta_3(y,u,r) |x_1|
        + \overline\mu_{2\psi}(\psi)\Big\{
        r^n\beta_4(x_1)[|\epsilon| + |\varpi|]
        + \beta_5(x_1,\hat\theta) |x_1| 
        + \theta \beta_6(x_1) |x_1|
        \Big\}
        + r^{n+1}\beta_7(x_1) [|\epsilon|^2+|\varpi|^2]
        + r^n\beta_8(x_1)\theta|\varpi||x_1|
        + \left|\frac{\partial\mu(t,x,\psi,u)}{\partial\psi} q_\psi(t,x,\psi,u,x_\Delta,\psi_\Delta,u_\Delta)\right|
        + \left|\frac{\partial\mu(t,x,\psi,u)}{\partial t} \right|$
        with
\preeqn\begin{align}
  \beta_1(y,u,r,\dot r,\hat\theta) &= \overline\mu_2(y,u)\{r\phi_{(1,2)}(x_1) + 
                                     r^n |K(x_1)|
                                     \nonumber\\
                                   &\hspace*{-0.22in}
                                     + \sum_{j=2}^{n-1}\phi_{(i,i+1)}(x_1) r^i + r^{n-1}n^{\frac{3}{2}}\Gamma(x_1)\}
                                     \nonumber\\
                                   &\hspace*{-0.22in}
                                     +nr^{n-1}\dot r |K(x_1)|
                                     + r^{n-1}\dot r |K(x_1)||D_c|
                                     \nonumber\\
                                   &\hspace*{-0.22in}
                                     + r^n\left|\frac{\partial K(x_1)}{\partial x_1}\right| |\hat\theta\vartheta_1(x_1) x_1|\phi_{(1,2)}(x_1)
                                     \nonumber\\
                                   &\hspace*{-0.22in}
  + r^{n+1}|K(x_1)||A_c(x_1)| 
  + r^{n+1}|K(x_1)||G(x_1)| 
                                     \nonumber\\
                                   &\hspace*{-0.22in}
  + r^n |K(x_1)|\hat\theta|\vartheta_1(x_1)+\vartheta_1'(x_1)x_1|\phi_{(1,2)}(x_1),
  \label{beta1defn}
\end{align}\posteqn
  $\beta_2(y,u,r,\hat\theta,\dot{\hat\theta}) = \overline\mu_2(y,u)\hat\theta\vartheta_1(x_1)\{\phi_{(1,2)}(x_1)+n\Gamma(x_1)\}  
  + r^{n-1}|K(x_1)| \hat\theta^2(|\vartheta_1(x_1)
  +\vartheta_1'(x_1)x_1|) \phi_{(1,2)}(x_1) |\vartheta_1(x_1)|
+r^{n-1}|K(x_1)||\dot{\hat\theta}\vartheta_1(x_1)|$,
  $\beta_3(y,u,r) = (n+1)\overline\mu_2(y,u)\Gamma(x_1) + r^n\frac{|K(x_1)| |G(x_1)|}{\phi_{(1,2)}(x_1)}\Gamma(x_1)         
  + r^{n-1}|K(x_1)|\hat\theta|\vartheta_1(x_1)+\vartheta_1'(x_1)x_1|\Gamma(x_1)$,
  $\beta_4(x_1) = \phi_{(1,2)}(x_1) + |K(x_1)| + \sum_{j=2}^{n-1}\phi_{(i,i+1)}(x_1) + n^{\frac{3}{2}}\Gamma(x_1)$,
  $\beta_5(x_1,\hat\theta) = \phi_{(1,2)}(x_1) \hat\theta \vartheta_1(x_1) + n\Gamma(x_1)\hat\theta\vartheta_1(x_1)$,
  $\beta_6(x_1) = (n+1)\Gamma(x_1)$,
  $\beta_7(x_1) = \frac{3}{2}\left|\frac{\partial K(x_1)}{\partial x_1}\right| \phi_{(1,2)}(x_1)$, and
  $\beta_8(x_1) = \left|\frac{\partial K(x_1)}{\partial x_1}\right| \Gamma(x_1)$.

Define
\preeqn\begin{align}
	 V_u &= \frac{1}{2r^n\Pi(r_u)}\ln(1+ u_d^2)
        \label{Vudefn}
\end{align}\posteqn
where $\ln$ denotes $log_e$; $\Pi$ is any function from ${\cal
  R}^+$ to ${\cal R}^+$ such that the following properties are satisfied:
$\Pi(a)\geq 1$ whenever $a\geq 1$; $\Pi(a)$ is a continuously
differentiable and monotonically increasing function over the interval
$[1,\infty)$;
a constant $\overline\Pi > 0$ exists such that
$\Pi(a)\leq \overline\Pi$ for all $a\in[1,\infty)$.
We utilize $\Pi'(r_u)$ to denote $\frac{\partial\Pi(r_u)}{\partial r_u}$.
An example of a function that satisfies the above properties
is $\Pi(a)=\tanh(k a) + 1$ with any constant $k > 0$. 
$V_u$ satisfies:
\preeqn\begin{align}
  \dot V_u 
&\leq
-r_u  \underline\mu\frac{u_d^2}{(1+u_d^2)r^n\Pi(r_u)}
\nonumber\\ &\quad
\!-\! \bigg(\frac{\Pi'(r_u)\dot r_u}{\Pi(r_u)} + \frac{n\dot r}{r}\bigg)\frac{\ln(1+u_d^2)}{2\Pi(r_u)r^n}
\nonumber\\ &\quad
+ k_{\psi 1}\frac{\alpha_\psi(\psi)}{r^{2n-\frac{3}{2}}}
+\frac{\Xi_{u1}(y,u,r,\dot r,\hat\theta,\dot{\hat\theta})}{\Pi(r_u)}\frac{u_d^2}{1+u_d^2} 
\nonumber\\ &\quad
+ r\tilde w_1(x_1)[|\epsilon|^2+|\varpi|^2]
\!+\! \frac{1}{r}\tilde q_1(x_1,\hat\theta) x_1^2 
              \!+\! \frac{1}{r}\theta^*\tilde q_2(x_1) x_1^2
              \nonumber\\ &\quad
                            + \frac{1}{r}\Big[c_3\theta^2  \!+\! c_2\hat\theta_\Delta^2\vartheta_1^2(x_{1,\Delta})\Big]\overline\mu_{1a}^2(x_{1,\Delta})x_{1,\Delta}^2
                            \nonumber\\ &\quad
                            + \frac{c_{\psi 1}}{2r^{2n-\frac{3}{2}}}\overline\mu_{1\psi}^2(\psi_\Delta)
                            \nonumber\\ &\quad
                            +\frac{c_1}{r^{2n-3}}\sum_{j=1}^nr_\Delta^{2j-2}(\varpi_{j,\Delta}^2+\epsilon_{j,\Delta}^2)\overline\mu_{1a}^2(x_{1,\Delta}) 
\label{Vudot2}
\end{align}\posteqnl
where
$r_\Delta$, $\hat\theta_\Delta$, $x_{1,\Delta}$, $\varpi_{j,\Delta}$, $\epsilon_{j,\Delta}$, and $\psi_\Delta$ denote the time-delayed versions (with time delay $\Delta$) of the corresponding signals (e.g., $x_{1,\Delta}(t)=x_1(t-\Delta)$, $\epsilon_{j,\Delta}(t)=\epsilon_j(t-\Delta)$). Also, 
$\tilde w_1(x_1) = c_1(1+\overline\mu_{1a}^2(x_1)) + \beta_7(x_1)  + \frac{\beta_4^2(x_1)}{c_{\psi 2}} + \frac{\beta_8^2(x_1)}{2c_4}$,
$\tilde q_1(x_1,\hat\theta) = c_2(1+\hat\theta^2\vartheta_1^2(x_1)\overline\mu_{1a}^2(x_1)) + \frac{\beta_5^2(x_1,\hat\theta)}{2c_{\psi 2}}$,
$\tilde q_2(x_1) = c_3(1+\overline\mu_{1a}^2(x_1)) + \frac{\beta_6^2(x_1)}{2c_{\psi 2}} + \frac{c_4}{2}$,
$k_{\psi 1} = 0.5(c_{\psi 1}+3c_{\psi 2})\overline k_\psi$,
and
\preeqnm\begin{align}
\Xi_{u1}(y,u,r,\dot r,\hat\theta,\dot{\hat\theta}) &=\frac{\overline\mu_2(y,u)}{r^n} +
\frac{1}{2c_{\psi_1}r^{\frac{3}{2}}} + \frac{\beta_1^2(y,u,r,\dot r, \hat\theta)}{2c_1r^{2n+1}}
\nonumber\\&\hspace*{-0.75in}
+ \frac{1}{c_{\psi 1}}\frac{\tilde{\overline\mu}_1^2(y,u)}{r^{\frac{3}{2}}}
+ \frac{1}{4c_2}\frac{\beta_2^2(y,u,r,\hat\theta,\dot{\hat\theta})}{r^{2n-1}}
+ \frac{1}{4c_3}\frac{\beta_3^2(y,u,r)}{r^{2n-1}}
\nonumber\\&\hspace*{-0.75in}
+\left[  \frac{1}{2c_3 r^{2n-1}} + \frac{n}{c_1r^3} + \frac{1}{2c_2 r^{2n-1}}\right]\tilde{\overline\mu}_1^2(y,u)
\end{align}\posteqn
where $c_1,\ldots,c_4$, $c_{\psi 1}$, and $c_{\psi 2}$ are any positive constants and  $\tilde{\overline\mu}_1(y,u)$ is defined as $[\overline\mu_1(y,u)+\tilde\mu_1(y,u)]$.
Since we have the upper bound
$|u_d| \leq \overline u_d(y,u,r,\varpi)$
where $\overline u_d(y,u,r,\varpi)$ is defined as
$\overline u_d(y,u,r,\varpi) \eqdef \Big[\overline\mu(y,u)+r^n|K(x_1)||\varpi|\Big]$,    
 the following inequality can be written for the term involving $u_d$ in \eqref{vxdot}:
   $- \frac{2c r\epsilon^T P_o B u_d}{r^{n-1}}  \leq
   r\lambda_{max}^2(P_o)  c^2 |\epsilon|^2
   + \frac{1}{r^{2n-3}}[1+\overline u_d^2(y,u,r,\varpi)]\ln(1+u_d^2)$. Note that $\overline u_d(y,u,r,\varpi)$ is a completely known function and involves only available variables.

\vspace*{0.03in}
\noindent{\bf Scaling of Lyapunov function for $\psi$ subsystem:}
Defining
    $\tilde V_\psi = \frac{V_\psi}{r^{2n-\frac{3}{2}}}$, we have
\preeqn\begin{align}
  \dot{\tilde V}_\psi
    &\leq
    -\frac{\alpha_\psi(|\psi|)}{r^{2n-\frac{3}{2}}} \!+\! \sum_{k=0}^1\Gamma_2(x_{1,k\Delta})\bigg\{\theta
    x_{1,k\Delta}^2 
    \nonumber\\&\quad
    + 3[|\varpi_{k\Delta}|^2+|\epsilon_{k\Delta}|^2]
    +\frac{3}{r_{k\Delta}^{\frac{1}{2}}}\hat\theta_{k\Delta}^2\vartheta_1^2(x_{1,k\Delta})x_{1,k\Delta}^2
    \nonumber\\&\quad
    +2\frac{\overline\gamma_s(x_{1,k\Delta})}{r_{k\Delta}^{2n-\frac{3}{2}}}u_{d,k\Delta}^2
    \nonumber\\&\quad
    +2r_{k\Delta}^{\frac{3}{2}}\overline\gamma_s(x_{1,k\Delta})|K(x_{1,k\Delta})|^2|\varpi_{k\Delta}|^2\bigg\}
\label{Vetadot1}
\end{align}\posteqnm

\noindent{\bf Overall composite Lyapunov function $V$ including terms to handle time-delayed terms in $\dot V_u$ and $\dot{\tilde V}_\psi$:}
Define
\preeqnm\begin{align}
          V &= V_{x} + c_\psi \tilde V_\psi + c_u V_u + \frac{1}{2c_\theta}(\hat\theta-\theta^*)^2
              \nonumber\\ &\quad
                            + \frac{1}{1-\overline\Delta}\int_{t-\Delta}^t\bigg\{
c_\psi\Gamma_2(x_1(\pi))\bigg\{\theta
    x_1^2(\pi) 
    \nonumber\\&\quad
    + 3[|\varpi(\pi)|^2+|\epsilon(\pi)|^2]
    +\frac{3}{r^{\frac{1}{2}}(\pi)}\hat\theta^2(\pi)\vartheta_1^2(x_1(\pi))x_1^2(\pi)
    \nonumber\\&\quad
                 +2r^{\frac{3}{2}}(\pi)\overline\gamma_s(x_1(\pi))|K(x_1(\pi))|^2|\varpi(\pi)|^2
    \nonumber\\&\quad
    +2\frac{\overline\gamma_s(x_1(\pi))}{r^{2n-\frac{3}{2}}(\pi)}u_d^2(\pi)
                 \bigg\}
                 \nonumber\\&\quad
                 +                            
\frac{c_u}{r(\pi)}\Big[c_3\theta^2  + c_2\hat\theta^2(\pi)\vartheta_1^2(x_1(\pi))\Big]\overline\mu_{1a}^2(x_1(\pi))x_1^2(\pi)
\nonumber\\&\quad
                                          +\frac{c_uc_1}{r^{2n-3}(\pi)}\sum_{j=1}^nr^{2j-2}(\pi)(\varpi_j^2(\pi)+\epsilon_j^2(\pi))\overline\mu_{1a}^2(x_1(\pi))
                            \nonumber\\ &\quad
                              + \frac{c_uc_{\psi 1}}{2r^{2n-\frac{3}{2}}(\pi)}\overline\mu_{1\psi}^2(\psi(\pi))
                                          \bigg\}
              d\pi
\end{align}\posteqnm
with $c_u$, $c_\theta$, and $c_\psi$ being any positive constants with furthermore $c_\psi$ additionally satisfying
  $c_\psi \geq c_u[2  k_{\psi 1}+\tilde\Delta c_{\psi 1}\overline k_\psi]$.
  Hence, using \eqref{Vudot2} and \eqref{Vetadot1},
we obtain
\preeqn\begin{align}
  \dot V 
&\leq
-r^2\frac{c\nu_o}{4}|\epsilon|^2
\!-\!r^2\frac{\nu_c}{4}\phi_{(2,3)}(x_1)|\varpi|^2
\!-\!\dot r c\ul\nu_o|\epsilon|^2
\!-\!\dot r \ul\nu_c|\varpi|^2
  \nonumber\\&\quad
-x_1\phi_{(1,2)}\vartheta \!-\! \frac{1}{2}\frac{\dot r}{r^2} x_1^2
-\frac{c_\psi \alpha_\psi(|\psi|)}{2 r^{2n-\frac{3}{2}}}
  \nonumber\\&\quad
-c_u r_u  \underline\mu\frac{u_d^2}{(1\!+\!u_d^2)r^n\Pi(r_u)}
- c_u \frac{\Pi'(r_u)\dot r_u}{2\Pi^2(r_u)r^n}\ln(1+u_d^2)
  \nonumber\\&\quad
-c_u \frac{n\dot r}{2r^{n+1}\Pi(r_u)}\ln(1+u_d^2)
  \nonumber\\&\quad
+\frac{1}{c_\theta}(\hat\theta-\theta^*)\dot{\hat\theta}
+ r\overline w_1(x_1,\hat\theta,\dot{\hat\theta})[|\epsilon|^2+|\varpi|^2]
  \nonumber\\&\quad
+ r^{\frac{3}{2}}\overline
w_2(x_1)[|\epsilon|^2+|\varpi|^2]
+\overline q_1(x_1)x_1^2 + \theta^* \overline q_2(x_1) x_1^2
\nonumber\\&\quad
+ \frac{\overline q_3(x_1,\hat\theta)\vartheta_1^2(x_1)x_1^2}{r^{\frac{1}{2}}}
+ \frac{\overline q_4(x_1,\hat\theta)}{r}x_1^2 + \theta^* \frac{\overline q_5(x_1)}{r}x_1^2
  \nonumber\\&\quad
+\frac{\overline \Xi_{u1}(y,u,r,\dot
  r,\hat\theta,\dot{\hat\theta})}{\Pi(r_u)}\frac{u_d^2}{1\!+\!u_d^2}
  \nonumber\\&\quad
+ \overline \Xi_{u2}(y,u,r,\varpi)\ln(1\!+\!u_d^2)
\end{align}\posteqn
where
  $\overline w_1(x_1,\hat\theta,\dot{\hat\theta}) = w_1(x_1,\hat\theta,\dot{\hat\theta}) 
  + \lambda_{max}^2(P_o)c^2
+  c_u [\tilde w_1(x_1)+c_1\tilde\Delta\overline\mu_{1a}^2(x_1)] +    3\tilde c_\psi \Gamma_2(x_1)$,
  $\overline w_2(x_1) = 2\tilde c_\psi \overline\gamma_s(x_1)\Gamma_2(x_1)|K(x_1)|^2$,
$\overline \Xi_{u1}(y,u,r,\dot
  r,\hat\theta,\dot{\hat\theta}) = c_u \Xi_{u1}(y,u,r,\dot
  r,\hat\theta,\dot{\hat\theta})$,
  $\overline \Xi_{u2}(y,u,r,\varpi) =
  \Big[\frac{1}{r^{2n-3}}
  + 2\tilde c_\psi \Gamma_2(x_1)\frac{\overline\gamma_s(x_1)}{r^{2n-\frac{3}{2}}}\Big]
  [1+\overline
  u_d^2(y,u,r,\varpi)]$,
$\overline q_1(x_1) = q_1(x_1)$, 
$\overline q_2(x_1) = q_2(x_1) + \tilde c_\psi \Gamma_2(x_1)$,
$\overline q_3(x_1,\hat\theta) = 3\tilde c_\psi\Gamma_2(x_1)\hat\theta^2$,
$\overline q_4(x_1,\hat\theta) = c_u [\tilde q_1(x_1,\hat\theta)+c_2\tilde\Delta\hat\theta^2\vartheta_1^2(x_1)\overline\mu_{1a}^2(x_1)]$, and
$\overline q_5(x_1) = c_u [\tilde q_2(x_1)+c_3\tilde\Delta\overline\mu_{1a}^2(x_1)]$, $\tilde c_\psi = c_\psi\frac{2-\overline\Delta}{1-\overline\Delta}$, and $\tilde\Delta=\frac{1}{1-\overline\Delta}$.

\begin{table*}[!t]
  \normalsize
\preeqnt\begin{align}
  R(x_1,\hat\theta,\dot{\hat\theta}) &= \max\bigg(\overline R,
  16\nu_a \overline w_1(x_1,\hat\theta,\dot{\hat\theta}),
  (16\nu_a\overline w_2(x_1))^2,
  \Big(\frac{4\overline
    q_3(x_1,\hat\theta)\vartheta_1(x_1)}{\hat\theta\phi_{(1,2)}(x_1)}\Big)^2,  
  \frac{4(\overline q_4(x_1,\hat\theta)+\hat\theta\overline q_5(x_1))}{\hat\theta\phi_{(1,2)}(x_1)\vartheta_1(x_1)}
    \bigg) \;\;\;\;\;
    \label{Rdefn}
\\
\Omega(x_1,\hat\theta,\dot{\hat\theta},r) &= \max\bigg(\overline\Omega,
2\nu_b r \overline w_1(x_1,\hat\theta,\dot{\hat\theta}),
2\nu_b r^{\frac{3}{2}}\overline w_2(x_1),
2\Big(\overline q_3(x_1,\hat\theta) \vartheta_1^2(x_1) r^{\frac{3}{2}} + r\overline
q_4(x_1,\hat\theta)+\hat\theta r \overline q_5(x_1)\Big)
    \bigg)
\label{Omegadefn}
\end{align}\posteqnm
\vspace*{-0.07in}
\hrule
\vspace*{-0.1in}
\end{table*}
\begin{table*}[!t]
\vspace*{-0.05in}
  \normalsize
\preeqn\begin{align}
  R_u(y,u,r,\dot r,\hat\theta,\dot{\hat\theta},\varpi) &= \max\bigg(\overline R_u,
  \frac{4r^n \overline \Xi_{u1}(y,u,r,\dot r,\hat\theta,\dot{\hat\theta})}{c_u\ul\mu},
  \frac{8r^n\overline\Pi (\overline \Xi_{u2}(y,u,r,\varpi)+\nu_u)(1+\overline u_d^2(y,u,r,\varpi))}{c_u\ul\mu}
    \bigg)
    \label{R2defn}
\\
\Omega_u(y,u,r,\dot r,\hat\theta,\dot{\hat\theta},\varpi,r_u) &= \max\bigg(\overline\Omega_u,
\frac{2\Pi^2(r_u)r^n}{c_u\Pi'(r_u)}\Big(\frac{\overline
  \Xi_{u1}(y,u,r,\dot
  r,\hat\theta,\dot{\hat\theta})}{\Pi(r_u)}+\overline \Xi_{u2}(y,u,r,\varpi)+\nu_u\Big)
    \bigg)
\label{Omega2defn}
       \end{align}\posteqnm
       \vspace*{-0.1in}
\hrule
\vspace*{-0.17in}
\end{table*}

\vspace*{0.03in}
\noindent{\bf Design of function $\vartheta_1$ appearing in \eqref{varthetadefn}:}
The function $\vartheta_1(x_1)$ is chosen as
  $\vartheta_1(x_1) = \frac{4}{\phi_{(1,2)}(x_1)}\Big(\frac{\overline
  q_1(x_1)+ \vartheta_1^*}{a_\theta}+\overline q_2(x_1) \Big)$
with $\vartheta_1^*$ being any positive constant.

\vspace*{0.03in}
\noindent{\bf Design of function $\overline Q_\theta$ appearing in dynamics of $\hat\theta$ in \eqref{thetahdot}:}
  $\overline Q_\theta(x_1) = c_\theta\Big(\overline
  q_2(x_1)+\frac{\overline q_5(x_1)}{r}\Big)x_1^2$.

\vspace*{0.03in}
\noindent{\bf Design of functions $\lambda$, $R$, $\Omega$, $R_u$, and $\Omega_u$ appearing in dynamics of scaling parameters in \eqref{rdot} and \eqref{r2dot}:}
 $\lambda:{\cal R}\rightarrow{\cal R}^+$ is chosen to be any
continuous function such that
$\lambda(s)=1$ for any $s>0$ and $\lambda(s)=0$ for
  any $s<-\epsilon_r$
with $\epsilon_r > 0$ being any constant. 
The functions $R$ and $\Omega$ are chosen as shown in \eqref{Rdefn} and \eqref{Omegadefn}, respectively,
with
$\overline R$ and $\overline\Omega$ being any nonnegative constants and
with
$\nu_a = \max\Big(\frac{1}{c\nu_o},\frac{1}{\nu_c\sigma}\Big)$
and
$\nu_b = \max\Big(\frac{1}{c\ul\nu_o},\frac{1}{\ul\nu_c}\Big)$.
 $R_u$ and $\Omega_u$ are chosen as shown in \eqref{R2defn} and \eqref{Omega2defn}, respectively,
where
$\overline R_u$ and $\overline\Omega_u$ are any nonnegative constants and
$\nu_u$ is any positive constant.

\vspace*{0.03in}
\noindent{\bf Analysis of closed-loop stability and asymptotic convergence:}
There are two possible cases: (i) $r \geq R(x_1,\hat\theta,\dot{\hat\theta})$; (ii) $r < R(x_1,\hat\theta,\dot{\hat\theta})$.
In case (i),  we have $r\geq
R(x_1,\hat\theta,\dot{\hat\theta})$, and in case (ii), we have $\dot r =
\Omega(x_1,\hat\theta,\dot{\hat\theta},r)$. In either of these cases, i.e., at all time instants $t$, we have
\preeqn\begin{align}
  \dot V 
&\leq
-r^2\frac{c\nu_o}{8}|\epsilon|^2
-r^2\frac{\nu_c}{8}\phi_{(2,3)}(x_1)|\varpi|^2
\!-\!\vartheta_1^* x_1^2 
\!-\!\frac{c_\psi \alpha_\psi(|\psi|)}{2 r^{2n-\frac{3}{2}}}
  \nonumber\\&\!\!\!
-c_u r_u  \underline\mu\frac{u_d^2}{(1+u_d^2)r^n\Pi(r_u)}
- c_u \frac{\Pi'(r_u)\dot r_u}{2\Pi^2(r_u)r^n}\ln(1+u_d^2)
  \nonumber\\&\!\!\!
+\frac{\overline \Xi_{u1}(y,u,r,\dot
  r,\hat\theta,\dot{\hat\theta})}{\Pi(r_u)}\frac{u_d^2}{1\!+\!u_d^2}
\!+\! \overline \Xi_{u2}(y,u,r,\varpi)\ln(1+u_d^2)
\label{eq:Vdot_a1}
\end{align}\posteqn
Now, similarly, there are two possible cases for $r_u$: (A) $r_u \geq R_u$; (B) $r_u < R_u$.
In case (A), we have $r_u\geq R_u$, and in case (B), we have $\dot r_u = \Omega_u$. 
It can be shown that
 in either of cases (A) and (B), i.e., at all time instants $t$,
\preeqn\begin{align}
   \dot V 
&\leq
-r^2\frac{c\nu_o}{8}|\epsilon|^2
-r^2\frac{\nu_c}{8}\phi_{(2,3)}(x_1)|\varpi|^2
-\vartheta_1^* x_1^2
\nonumber\\&\quad
-\frac{c_\psi \ul V_\psi }{2} \tilde V_\psi(\psi)
-\nu_u \ln(1+u_d^2).
\label{vdot3_2}
       \end{align}\posteqnmm
       From \eqref{vdot3_2} and the system dynamics, the closed-loop stability and asymptotic convergence properties can be inferred.
Firstly, to show that solutions exist for all time, 
consider 
the maximal interval of existence of solutions of the closed-loop dynamic system to be $[0,t_f)$ with some $t_f > 0$.
Then, from (\ref{vdot3_2}), it can be seen that
$V$ is bounded
on $[0,t_f)$ and that therefore,  
from the definition
of $V$ and the dynamics of the closed-loop system, a process of signal chasing can be used to show that all closed-loop signals remain bounded on $[0,t_f)$. %
Therefore,
  solutions
exist for all time, i.e., $t_f=\infty$. 
Also, from (\ref{vdot3_2}),
$\epsilon$, $\varpi$, $x_1$, $\psi$, and $u_d$, and therefore,
$x_1,\ldots,x_n,\psi,\hat x_2,\ldots, \hat x_n$ 
go to
zero
asymptotically
as $t\rightarrow \infty$.  %

\section{Conclusion}
We considered a general class of uncertain nonlinear systems with input unmodeled dynamics with time-varying time delays and showed that
a dynamic scaling based robust adaptive output-feedback controller
can be designed to globally stabilize this uncertain nonlinear system.
While the input unmodeled dynamics subsystem involves uncertain (and time-varying) time delays both on the state and control input signals entering into this subsystem, 
the developed control design methodology is itself delay-independent in two important senses: firstly, the controller does not utilize any delayed versions of the measured output, input, or controller internal variable signals; secondly, the control design does not require knowledge of the actual time delay magnitudes (specifically, only requires an upper bound on rate of change of the time delay magnitude).
In further research, applicability of the proposed techniques to more general classes of nonlinear systems (e.g., nontriangular systems, more general structures of appended dynamics, etc.) are being considered.

\renewcommand{\baselinestretch}{0.98}
\bibliographystyle{IEEEtran}
\bibliography{refs}

\begin{thebibliography}{00}
\providecommand{\url}[1]{#1}
\csname url@rmstyle\endcsname
\providecommand{\newblock}{\relax}
\providecommand{\bibinfo}[2]{#2}
\providecommand\BIBentrySTDinterwordspacing{\spaceskip=0pt\relax}
\providecommand\BIBentryALTinterwordstretchfactor{4}
\providecommand\BIBentryALTinterwordspacing{\spaceskip=\fontdimen2\font plus
\BIBentryALTinterwordstretchfactor\fontdimen3\font minus
  \fontdimen4\font\relax}
\providecommand\BIBforeignlanguage[2]{{%
\expandafter\ifx\csname l@#1\endcsname\relax
\typeout{** WARNING: IEEEtran.bst: No hyphenation pattern has been}%
\typeout{** loaded for the language `#1'. Using the pattern for}%
\typeout{** the default language instead.}%
\else
\language=\csname l@#1\endcsname
\fi
#2}}
\providecommand{\BIBdecl}{\relax}
\BIBdecl


 \bibitem{KKK95}
 M.~Krsti\'c, I.~Kanellakopoulos, and P.~V. Kokotovi\'c, \emph{Nonlinear and
   Adaptive Control Design}.\hskip 1em plus 0.5em minus 0.4em\relax New York:
   Wiley, 1995.

 \bibitem{Isi99}
 A.~Isidori, \emph{Nonlinear Control Systems II}. London, UK: Springer, 1999. 

 \bibitem{Kha01}
 H.~Khalil, \emph{Nonlinear Systems}. Upper Saddle River, NJ: Prentice Hall, 2001. 


  \bibitem{PW96}
  L. Praly and Y. Wang, ``Stabilization in spite of matched unmodeled dynamics and an equivalent definition of input-to-state stability,'' \emph{Math. of Control, Signals and Systems}, vol. 9, no. 1, pp. 1-33, 1996.

\bibitem{SJK97b}
R.~Sepulchre, M.~Jankovi{\'c}, and P.~Kokotovi{\'c}, \emph{Constructive
  Nonlinear Control}.\hskip 1em plus 0.5em minus 0.4em\relax London, UK:
  Springer Verlag, 1997.
  

\bibitem{RGLS08}
B.~Ren, S.~S.~Ge, T.~H.~Lee, and C.-Y.~Su, ``Adaptive neural control for uncertain nonlinear systems in pure-feedback form with hysteresis input,'' in \emph{Proc. of
  the IEEE Conf. on Decision and Control}, Cancun, Mexico, Dec. 2008, pp.  86--91.

\bibitem{Jan05}
M.~Jankovic, ``Stabilization of nonlinear time delay systems with
  delay-independent feedback,'' in \emph{Proc. of the American Control
  Conference}, Portland, OR, June 2005, pp. 4253--4258.

\bibitem{MB06}
F.~Mazenc and P.-A. Bliman, ``Backstepping design for time-delay nonlinear
  systems,'' \emph{IEEE Trans. on Automatic Control}, vol.~51, no.~1, pp.
  149--154, Jan. 2006.

\bibitem{MML07}
F.~Mazenc, M.~Malisoff, and Z.~Lin, ``On input-to-state stability for nonlinear
  systems with delayed feedbacks,'' in \emph{Proc. of the American
  Control Conference}, New York, NY, July 2007, pp. 4804--4809.





\bibitem{IR03}
A.~Ilchmann and E.~P. Ryan, ``On gain adaptation in adaptive control,''
  \emph{IEEE Trans. on Automatic Control}, vol.~48, no.~5, pp. 895--899,
  May 2003.

\bibitem{Pra03}
L.~Praly, ``Asymptotic stabilization via output feedback for lower triangular
  systems with output dependent incremental rate,'' \emph{IEEE Trans. on
  Automatic Control}, vol.~48, no.~6, pp. 1103--1108, June 2003.


\bibitem{KK01a} 
  P.~Krishnamurthy and F.~Khorrami, ``Global adaptive output feedback tracking for nonlinear systems linear in unmeasured states,'' in \emph{Proc. of the American Control Conference}, Arlington, VA, June 2001, pp.~4814--4819.

\bibitem{KK01b} 
  P.~Krishnamurthy and F.~Khorrami, ``Decentralized control of large-scale nonlinear systems in generalized output-feedback canonical form,'' in \emph{Proc. of the IEEE Conference on Decision and Control}, Orlando, FL, Dec. 2001, pp.~1322--1327.

  
\bibitem{KKC02}
P.~Krishnamurthy, F.~Khorrami, and R.~S. Chandra, ``Global high-gain-based
  observer and backstepping controller for generalized output-feedback
  canonical form,'' \emph{IEEE Trans. on Automatic Control}, vol.~48,
  no.~12, pp. 2277--2284, Dec. 2003.

\bibitem{KK03b}
   P.~Krishnamurthy and F.~Khorrami, ``Decentralized control and disturbance attenuation for large-scale nonlinear systems in generalized output-feedback canonical form,'' \emph{Automatica}, vol.~39, no.~11, pp.~1923--1933, Nov. 2003.

  
  \bibitem{KK03}
P.~Krishnamurthy and F.~Khorrami, ``Robust adaptive control for non-linear systems in generalized output-feedback canonical form,'' \emph{International Journal of Adaptive Control and Signal Processing}, vol.~17,
no.~4, pp.~285--311, May 2003. 


\bibitem{Ito06}
H.~Ito, ``State-dependent scaling problems and stability of interconnected iISS and ISS systems,'' \emph{IEEE Trans. on
  Automatic Control}, vol.~51, no.~10, pp. 1626--1643, Oct. 2006.


\bibitem{Lin09}
Z.~Lin, ``Low gain and low-and-high gain feedback: A review and some recent results,'' in \emph{Proc. of
  the IEEE Conf. on Decision and Control}, Shanghai, China, Dec. 2009, pp.
  lii--lxi.

  \bibitem{KK09a}
   P.~Krishnamurthy and F.~Khorrami, ``Application of a dynamic high-gain scaling methodology to servocompensator design,'' \emph{International Journal of Robust and Nonlinear Control}, vol.~19, no.~8, pp.~937--964, May 2009.


\bibitem{NK11}
  S. Nazrulla and H. K. Khalil, ``Robust Stabilization of Non-Minimum Phase Nonlinear Systems Using Extended High-Gain Observers,'' \emph{IEEE Transactions on Automatic Control}, vol.~56, no.~4, pp.~802--813, April 2011.


\bibitem{KK03a} 
  P.~Krishnamurthy and F.~Khorrami, ``A dual high gain controller for the uncertain generalized output-feedback canonical form with appended dynamics driven by all states,'' in \emph{Proc. of the American Control Conference}, Denver, CO, June 2003, pp. 4766--4771.
  
\bibitem{KK04e}
P.~Krishnamurthy and F.~Khorrami, ``A high-gain scaling technique for adaptive
  output feedback control of feedforward systems,'' \emph{IEEE Trans. on
  Automatic Control}, vol.~49, no.~12, pp. 2286--2292, Dec. 2004.

\bibitem{KK04f}
P.~Krishnamurthy and F.~Khorrami, ``Dynamic high-gain scaling: state and output feedback with application
  to systems with {ISS} appended dynamics driven by all states,'' \emph{IEEE
  Trans. on Automatic Control}, vol.~49, no.~12, pp. 2219--2239, Dec.
  2004.


\bibitem{KAP06}
G.~Kaliora, A.~Astolfi, and L.~Praly, ``Norm estimators and global output
  feedback stabilization of nonlinear systems with {ISS} inverse dynamics,''
  \emph{IEEE Trans. on Automatic Control}, vol.~51, no.~3, pp. 493--498,
  March 2006.

 \bibitem{KK07f}
 P.~Krishnamurthy and F.~Khorrami, ``High-gain output-feedback control for
   nonlinear systems based on multiple time scaling,'' \emph{Systems and Control
   Letters}, vol.~56, no.~1, pp. 7--15, Jan. 2007.


\bibitem{KK07a}
P.~Krishnamurthy and F.~Khorrami, ``Generalized state scaling and applications to feedback, feedforward,
  and non-triangular nonlinear systems,'' \emph{IEEE Trans. on Automatic
  Control}, vol.~52, no.~1, pp. 102--108, Jan. 2007.

\bibitem{KK07b} 
  P.~Krishnamurthy and F.~Khorrami, ``Adaptive output-feedback control of feedforward systems with uncertain parameters coupled with all states,'' in \emph{Proc. of the American Control Conference}, New York, NY, July 2007, pp. 480--485.

\bibitem{KK07h}
P.~Krishnamurthy and F.~Khorrami, ``Feedforward systems with {ISS} appended dynamics: Adaptive
  output-feedback stabilization and disturbance attenuation,'' \emph{IEEE
  Trans. on Automatic Control}, vol.~53, no.~1, pp. 405--412, Feb. 2008.

\bibitem{KK08a}
  P.~Krishnamurthy and F.~Khorrami, ``Dual high-gain-based adaptive output-feedback control for a class of nonlinear
systems,'' International Journal of Adaptive Control and Signal Processing, vol. 22,
no. 1, pp. 23--42, Feb. 2008.

\bibitem{KK10b}
P.~Krishnamurthy and F.~Khorrami, ``Output-feedback control of nonlinear delayed systems: A dynamic high-gain scaling approach,'' Dynamics of Continuous, Discrete and Impulsive Systems. Special Issue. vol. 17, no. 6, pp. 909-934, 2010. 


\bibitem{KK11b}
P. Krishnamurthy and F. Khorrami, ``Adaptive output-feedback control of a general class of uncertain feedforward systems via a dynamic scaling approach,'' \emph{IET Control Theory and Applications}, vol.~5, no.~5, pp.~681--692, March 2011. 

  \bibitem{KK13a} 
P.~Krishnamurthy and F.~Khorrami, ``A singular perturbation based global dynamic high gain scaling control design for systems with nonlinear input uncertainties,'' \emph{IEEE Trans. on Automatic Control}, vol.~58, no.~10, pp.~2686--2692, Oct. 2013. 


    \bibitem{KK14c}
P.~Krishnamurthy and F.~Khorrami, ``Dynamic High Gain Scaling Based Output Feedback for Nonlinear Systems with Time-Delayed Input Unmodeled Dynamics,'' in \emph{Proc. of the European Control Conference}, Strasbourg, France, June 2014, pp. 1777--1782.

  \bibitem{KK14b} 
P.~Krishnamurthy and F.~Khorrami, ``A dynamic scaling based control redesign procedure for uncertain nonlinear systems with input unmodeled dynamics,'' in \emph{Proc. of the IFAC World Congress}, Cape Town, South Africa, Aug. 2014. 

\bibitem{KK15a}
P.~Krishnamurthy and F.~Khorrami, ``Multiple-dynamic-scaling output-feedback control for uncertain
  strict-feedback-like systems with input unmodeled dynamics,'' in \emph{Proc.
  of the American Control Conference}, Chicago, IL, July 2015, pp. 2685--2690.

\bibitem{SOK84}
V.~R. Saxena, J.~O'Reilly, and P.~V. Kokotovic, ``Singular perturbations and
  time-scale methods in control theory: survey 1976-1983,'' \emph{Automatica},
  vol.~20, no.~3, pp. 273--293, 1984.

\bibitem{HLS07}
N.~Hovakimyan, E.~Lavretsky, and A.~Sasane, ``Dynamic inversion for
  nonaffine-in-control systems via time-scale separation: Part I,''
  \emph{Journal of Dyn. and Control Systems}, vol.~13, no.~4, pp.
  451--465, Oct. 2007.

\bibitem{CA09}
A.~Chakrabortty and M.~Arcak, ``Time-scale separation redesigns for
  stabilization and performance recovery of uncertain nonlinear systems,''
  \emph{Automatica}, vol.~45, no.~1, pp. 34--44, 2009.

\bibitem{KK04g} 
  P.~Krishnamurthy and F.~Khorrami, ``Conditions for uniform solvability of parameter-dependent Lyapunov equations with applications,'' in \emph{Proc. of the American Control Conference}, Boston, MA, July 2004, pp. 3896-3901.
  
\bibitem{KK06g}
P.~Krishnamurthy and F.~Khorrami, ``On uniform solvability of
  parameter-dependent Lyapunov inequalities and applications to various
  problems,'' \emph{SIAM Journal on Control and Optimization}, vol.~45, no.~4, pp.~1147-1164, Sep. 2006.



  
\end{thebibliography}

\end{document}